\title{Computing ``Small" $1$--Homological Models\\
for Commutative Differential Graded Algebras\thanks{
Proceedings of the 3rd Workshop on Computer Algebra and Scientific Computing CASC'00,
Springer--Verlag (2000), 87--100.}}
\author{\normalsize CHATA group\thanks{\'{A}lvarez V.,  Armario J.A., Frau M.D.,
Gonz\'{a}lez--D\'{\i}az R., Jim\'{e}nez M.J., Real P. and Silva B. Authors are partially supported by the PAICYT research project
FQM--0143 from Junta de Andaluc\'{\i}a and the DGES--SEUID research
projects
PB98--1621--C02--02
from Education and Science Ministry (Spain).}\\
\normalsize Dept. of Applied Math.,\\ \normalsize University  of Seville, Spain,\\
\normalsize  real@us.es}
 \date{}

\documentstyle[emlines2,bezier,12pt]{article}
\newcommand{\N}{\bf N}

\newcommand{\scst}{\scriptscriptstyle}
\newcommand{\ot}{\otimes}
\newcommand{\ra}{\rightarrow}
\newtheorem{procedure}{Procedure}{\bfseries}{\itshape}
\newtheorem{prp}{Proposition}[section]
\newtheorem{lemma}[prp]{Lemma}
\newtheorem{theorem}[prp]{Theorem}
\newtheorem{definition}[prp]{Definition}

\begin{document}
\sloppy
\maketitle
\begin{abstract}
We use homological perturbation machinery specific
for the algebra category \cite{Rea96} to give
an  algorithm for computing the differential structure of a small
$1$--homological model for commutative differential graded algebras
(briefly, CDGAs). The
complexity of the procedure is studied and a computer package in Mathematica is described
for determining such models.
\end{abstract}

\section{Introduction}

The description of efficient algorithms for
homological computation can be
considered to be a very important topic in Homological Algebra.
These algorithms can be  used mainly in the resolution of problems in
Algebraic Topology; but this subject also impinges directly
on the development of diverse areas such as Combinatorial Designs, Code Theory, Concurrency
Theory or Cohomological Physics.

Starting from a finite CDGA $A$, we establish an algorithm for
obtaining an ``economical'' $1$--homological model  $hBA$,
in the sense
that the number of algebra generators of $hBA$ is less than that of
the reduced bar construction $\bar{B}(A)$.
In order to get the $1$--homology of $A$, we would need to compute
 the homology groups of the model $hBA$. This computation can be reduced to
a simple problem of Linear Algebra (see \cite{Mun84} for a
complete explanation of this method).

Our main technique is  homological
  perturbation machinery \cite{GL89,GLS91,HK91}. Homological
  Perturbation Theory is often used to replace given chain
complexes by other smaller, homotopic
chain complexes which are more readily computable. An essential notion in this theory is that of
contraction. A contraction $c=(f,g,\phi)$ between two differential
graded modules
$(N,d_{\scst N})$ and
$(M,d_{\scst M})$ is a special homotopy equivalence between both modules
such that
 the corresponding homology groups are isomorphic.
 The morphisms $f$, $g$ and $\phi$ are
 called {\em
 projection}, {\em inclusion} and
 {\em homotopy} of the contraction, respectively. The Basic
Perturbation
Lemma  is
the heart of this theory and states
that given a contraction $c= (f,g,\phi)$, and a perturbation
$\delta$
of $d_{\scst N}$ (that is, $(d_{\scst N}+\delta)^{2}=0$),
then there exists
a new contraction $c_{\delta}=(f_{\delta}, g_{\delta}, \phi_{\delta})$
from $(N,d_{\scst N}+ \delta)$
to $(M,d_{\scst M}+d_{\delta})$, satisfying
\begin{eqnarray*}
f_{\delta}  =  f
(1- \delta \Sigma_{c}^{\delta} \phi)\,,\quad
g_{\delta}  =  \Sigma_{c}^{\delta} g\,,\quad
\phi_{\delta}  =  \Sigma_{c}^{\delta} \phi\end{eqnarray*}
\begin{eqnarray}
\label{dedelta}
d_{\delta} =  f \delta  \Sigma_{c}^{\delta}  g\,,
\end{eqnarray}

\noindent where
$\Sigma_{c}^{\delta}\;  =\;
\displaystyle\sum_{i\geq 0} (-1)^i \; (\phi\delta)^i \; =\;
1- \phi \delta + \phi\delta\phi\delta - \cdots +
(-1)^i (\phi\delta)^i + \cdots$.

It is necessary to emphasize that a nilpotent condition for the
composition $\delta \phi$ is required for guaranteeing  the
finiteness of the formulas.

The basic idea we use in this paper is the establishment (via
composition, tensor product or perturbation of contractions) of an
explicit contraction from an initial differential graded module $N$ to a free
differential graded module $M$ of finite type, so that the homology of $N$
is computable from that of $M$.

This ``modus operandi'' has been used by the authors in previous works
\cite{AAGR98,GR99,ARS99}.

Working in the context of CDGAs, Homological Perturbation Theory
 immediately supplies  a general
algorithm computing the $1$--homology of these objects at graded
module level. Nevertheless, this procedure,
already presented by Lambe in \cite{Lam92}, bears, in general,  high computational
charges and actually restricts its application to the low
dimensional homological calculus.

   This algorithm is refined, taking advantage
   of the multiplicative structures, in  \cite{AARS97}. More precisely,
   the Semifull Algebra Perturbation
Lemma \cite[Sec. 4]{Rea96} is used for designing the algorithm {\bf Alg1}. The input of this
method is a
CDGA $A$ given in the form of a ``twisted'' tensor product of $n$
exterior and polynomial
algebras, and the output  is a contraction $c_{\delta}$
(produced via  perturbation) from the reduced bar
construction $\bar{B}(A)$  to a smaller
differential graded algebra $hBA$, which is free and of finite type
as a graded module. In this case, we say that
the pair $\{c_{\delta},hBA\}$ (or, simply, $hBA$)
is a $1$--homological model of $A$.
Taking advantage
of the fact that the differential $d_{\scst hBA}$ of $hBA$ is
a derivation, that is, a morphism compatible with the product
of the $1$--homological model, it is only necessary to know the value of this morphism
applied to the generators of the model (let us observe that there are $n$
algebra generators). This implies a substantial improvement in the computation
of the differential on the small model $hBA$.

     In spite of this improvement, the computational cost for determining the
     morphism $d_{\scst hBA}$ applied to an algebra generator of $hBA$ is
     enormous, since the differential $d_{\scst hBA}$ follows the
     formula  (\ref{dedelta}) and the homotopy $\phi$ of $c$
     has an essentially exponential nature not only in time but
     also in space.

   We develop  some techniques, which comprise what we call
   {\em Inversion Theory} and which first appears in \cite{Rea96}.
   In consequence, we refine
the formula for $\phi$,
which is involved in the description of $d_{\scst hBA}$.
This study is based on the observation that the
the projection
 $f$ applied to certain elements (those ``with inversions") is always
 null. It follows that a not insignificant
number of terms in the formula of the
 morphism $\phi$ can be eliminated in the composition
 $f \delta (\phi\delta)^i g$, which appears  in the formula (\ref{dedelta}) of
 $d_{\delta}$. In such a way, we
 derive  an upgraded algorithm {\bf Alg2}.

The article is organized as follows: Notation and terminology
are introduced in Section \ref{preliminares}. In Section
\ref{algoritmo} the algorithm {\bf Alg1}, which was described in \cite{AARS97} is
recalled. Our contribution starts in Section \ref{inversiones}
which is devoted to explaining Inversion Theory and describing the algorithm {\bf
Alg2}.
An analysis of the complexity of {\bf Alg1} and {\bf Alg2}
 for computing the
differential structure of the small $1$--homological model $hBA$
is carried out in Section \ref{complejidad} and
a comparison between both algorithms is given.
Finally, in Section \ref{ejemplos}
we also give several examples illustrating the implementation of
{\bf Alg2} carried out using Mathematica 3.0.

\section{Preliminaries}\label{preliminares}

Although relevant notions of Homological Algebra are explained through the
exposition of this paper, most  common concepts are not explicitly given
(they can be found, for instance, in \cite{Mac95} or \cite{Wei94}).

Let $\Lambda$ be a commutative ring with the non zero unit, which will
be considered to be the ground ring. A {\em DGA--module}
$(M,d_{\scst M},\xi_{\scst M},\eta_{\scst M})$ is a
module endowed with:
\begin{itemize}
\item  A graduation, that is,
$M=\displaystyle\oplus_{n\in \N}M_n$.
\item A differential, $d_{\scst M}:M\rightarrow M$, which decreases the degree
by one and satisfies $d_{\scst M}^2=0$.
\item An augmentation,
$ \xi_{\scst M}:M\rightarrow\Lambda$, with $\xi_{\scst M}d_1=0$.
\item A coaugmentation, $\eta:\Lambda\rightarrow M$, with
$\xi_{\scst M}\eta_{\scst M}=1_{\scst \Lambda}$.
\end{itemize}
We will respect
Koszul conventions.
The {\em homology} of a differential graded module M,
is a graded module $H_*(M)$, where $H_n(M)=$Ker$\,d_n/$Im$\,d_{n+1}$.
We
are specially interested on CDGAs,
$(A,d_{\scst A},*_{\scst A},\xi_{\scst A},\eta_{\scst A})$
which are  differential graded modules endowed with a
product, $*_{\scst A}$, that is commutative in a graded sense.
A morphism  $\delta:A\rightarrow A$ which decreases the degree by one, is
a {\em derivation} if $\delta *_{\scst A}=*_{\scst A}(1\ot\delta+\delta\ot
1)$.

Three particular algebras are of special interest in the development of
this
paper: exterior, polynomial and divided power algebras. Let
$n$ be a fixed non--negative integer.
\begin{itemize}
\item The {\em exterior algebra}
$E(x,2n+1)$ is the graded
algebra with generators $1$ and $x$ of degrees $0$ and $2n+1$,
respectively, and the trivial product, that is, $x\cdot x=0$ and $x \cdot 1=x$.
\item The {\em polynomial algebra}
$P(y,2n)$ consists in the graded algebra with generators $1$ of
degree $0$ and  $y$ of degree $2n$. The
product is the usual one in polynomials, i.e.: $y^i \cdot
y^j=y^{i+j}$, for non negative integers $i$ and $j$.
\item Finally, the {\em divided
power algebra} $\Gamma (y,2n)$ is the graded algebra with generators $1$ and
$y$ ( $y=y^{(1)}$) of respective degrees $0$ and $2n$. The product is defined by the rules
$\mbox{$y^{(i)} \cdot y^{(j)} =
\left(\begin{array}{c} i+j \\ i\end{array} \right) \; y^{(i+j)}$}$,
 $i$ and $j$ being non--negative integers.
 \end{itemize}
Each one of these three types of
algebras can be considered as a CDGA with the trivial differential.

Now, we shall recall a standard
algebraic tool which allows us to preserve
the product structure of the initial
CDGA through the procedure of homological computation.
The {\it reduced bar construction} \cite{Mac95} associated to a CDGA $A$ is
defined as the differential graded module $\bar{B}(A)$:

$$\bar{B}(A)= \Lambda \oplus \mbox{Ker }\xi_{\scst A} \oplus (\mbox{Ker }\xi_{\scst A}
\otimes \mbox{Ker }\xi_{\scst A} ) \oplus \cdots \oplus (\mbox{Ker }\xi_{\scst A}
\otimes \cdots \otimes
\mbox{Ker }\xi_{\scst A} ) \oplus \cdots \,.$$

An element from $\bar{B}(A)$ is denoted by
$\bar{a}=[a_1 \vert \cdots \vert a_n]$.
There is a {\em tensor graduation} $|\;\;|_t$
 given by $|[a_1 | \cdots | a_n]|_t= \sum_{i=1}^{n}|a_i|$, as well
 as a {\em simplicial graduation} $\vert \;\;
\vert_s $, which is defined by
$\vert \bar{a} \vert_s= \vert [a_1\vert \cdots \vert a_n] \,
\vert_s = n.$
 The total degree of $\bar{a}$ is given by $\vert
\bar{a}\vert = \vert \bar{a}\vert_t +\vert \bar{a} \vert_s $.

The total differential is
given by the sum of the tensor and simplicial differentials.
The {\em tensor differential } is
defined by:
$$d_t[a_1| \cdots |a_n]=\displaystyle -\sum
_{i=1}^n(-1)^{|[a_1|\cdots | a_{i-1}]|}[a_1| \cdots |
d_{\scst
A}a_i | \cdots | a_n]\,.$$
The {\em simplicial differential} acts by cutting down the simplicial degree
by using the product given in $A$.

When the algebra $A$ is commutative, it is possible to define a
multiplicative structure on $\bar{B}(A)$ (via an operator called
the {\it  shuffle product}), so that the reduced bar construction
also becomes a CDGA.

Given two non--negative integers $p$ and $q$ , a
{\em $(p,q)$--shuffle} is defined as a
permutation $\pi$ of the set  $\{ 0,
\ldots ,p+q-1 \}$, such that $\pi (i) < \pi (j)$ when $0\leq i<j \leq p-1$
or $p \leq
i<j \leq p+q-1$ is the case.

Let us observe that there are $\left(\begin{array}{c} p+q \\
p\end{array} \right)$
different $(p,q)$--shuffles.

So, given a CDGA $A$, the
{\em shuffle product } $\star: \bar{B}(A) \ot \bar{B}(A) \longrightarrow
\bar{B}(A)$, is defined by:
$$[a_1| \cdots |a_p] \star [b_1| \cdots |b_q]= \displaystyle
\sum _{\scst \pi\in \{(p,q)--shuffles\}}
(-1)^ {\varepsilon ( \pi ,a,b)} [c_{\pi(0)}| \cdots |c_{\pi(p+q-1)}]\,;$$
where $(c_{0}, \ldots ,c_{p-1},c_{p}, \ldots ,c_{p+q-1})
 =
 (a_1, \ldots ,a_p,b_1,
\ldots b_q)$ and
$$\varepsilon ( \pi ,a,b)= \sum _{\pi (i)> \pi (p+j)}|[a_i]|
|[b_j]|\,.$$

Let $n$ be a non--negative integer.
The $n$--homology  of a CDGA $A$ (see \cite{Mac95})
consists in the homology groups
of the iterated
reduced bar construction
$\bar{B}^{n} (A)=
\bar{B}(\bar{B}^{n-1} (A) ) $, being   $\bar{B}^0(A)=A$.

Let $\{ A_i\}_{i \in I}$ be a set of
 CDGAs.
A  {\em twisted tensor product}
$\tilde{\ot}^{\rho}_{i \in I} A_i$ is a CDGA satisfying the following
conditions:
\begin{enumerate}
\item [\bf i)] $\tilde{\ot}^{\rho}_{i \in I} A_i$ coincides with the tensor
product $\ot_{i \in I} A_i$ as a graded algebra.
\item [\bf ii)] The differential operator consists in the sum of the
differential of the banal tensor product and a derivation
$\rho$.\end{enumerate}

A {\em contraction} $c:\{N,M,f,g,\phi\}$ \cite{EM53,HK91}, also denoted by
$(f, g, \phi ):N \stackrel
{c}{\Rightarrow} M$, from a differential graded module
$(N, d_{\scst N})$ to a differential graded module $(M, d_{\scst M})$ consists
in a homotopy
equivalence determined by three morphisms $f$, $g$ and $\phi$;
 $f:N_*
\rightarrow M_*$ (projection) and $g:M_* \rightarrow N_*$ (inclusion) being two
differential graded module morphisms and $\phi: N_* \rightarrow N_{*+1}$ a homotopy
operator. Moreover, these data are required to satisfy the following rules:
$$fg=1_{\scst M}\,, \quad
\phi d_{\scst N}+d_{\scst N}\phi + gf=1_{\scst
N}\,,\quad f \phi =0\,, \quad \phi g=0\,, \quad\quad \phi \phi =0\,.$$

 There are two basic
operations between contractions which give place to new
contractions:
tensor product and composition of contractions.

In this paper we use a  particular type of contraction between
CDGAs. Given two CDGAs $A$ and $A'$,
a {\it semifull algebra contraction}
$(f,g,\phi):A\Rightarrow A'$
\cite{Rea96} consists of an inclusion $g$
that is a
morphism of CDGAs, a quasi--algebra projection $f$
and a quasi--algebra
homotopy $\phi$.
We recall that
\begin{enumerate}
\item The projection $f$ is said to be a {\it quasi--algebra projection}
whenever the following conditions hold:
\begin{equation} \label{qapro}
f ( \phi *_{\scst A} \phi )=0, \quad f ( \phi *_{\scst A}
g)=0, \quad f (g *_{\scst A} \phi )=0\,.
\end{equation}
\item The homotopy operator $\phi$ is said to be a {\it quasi--algebra
homotopy} if
\begin{equation} \label{qahom}
\phi ( \phi *_{\scst A} \phi )=0, \quad \phi ( \phi
*_{\scst A} g)=0, \quad \phi (g *_{\scst A} \phi )=0\,.
\end{equation}
\end{enumerate}

The class
of all semifull algebra contractions is closed under composition and tensor
product of contractions. Moreover, this class is closed under
perturbation.

\begin{theorem}\label{lph}\cite{Rea96}

Let $c:\{N, M, f, g, \phi\}$ be a semifull algebra contraction and $\delta:
N \ra N$ be a perturbation--derivation of $d_{\scst N}$.
Then, the perturbed contraction $c_{\delta}$, is a
new semifull algebra contraction.
\end{theorem}

\section{Computability of the $1$--Homology of CDGAs. First
Algorithm}\label{algoritmo}

Here we recall the algorithm described in \cite{AARS97} for the
computation of a $1$--homological model of a CDGA.

It is commonly known that every CDGA $A$
``factors'',
up to homotopy equivalence,
into a tensor product of exterior and polynomial
algebras endowed with a
differential--derivation; in the sense that there exists a homomorphism
connecting both structures, which induces an isomorphism in
homology.

In fact, our input is  a twisted tensor product of
algebras $A=\tilde{\ot}^{\rho}_{i\in I} A_i$ where $I$ denotes a
finite set of indices, $\rho$  is a differential--derivation and $A_i$  an exterior
or a polynomial algebra, for every $i$.
In our algorithmic approach,  we encode $A$ by
\begin{enumerate}
\item a  sequence of non--negative integers
$n_1\leq n_2\leq\cdots\leq n_k$, such that $n_i$ represents the degree
of the algebra generator $x_i$ of $A_i$;
\item a $k$--vector $\bar{v}=(v_1,v_2,\dots,v_k)$, such that
$v_i$ is $\rho(x_i)$ for all $i$.
\end{enumerate}

The principal
goal is to obtain a ``chain''
of semifull algebra contractions starting at the
reduced bar construction $\bar{B}(A)$ and ending up at
a smaller free (as a module) CDGA. In that way, we determine
a $1$--homological model for $A$.

Now, we consider the following three
{\bf semifull} algebra contractions
which are  used, firstly,
to find the structure of a graded module of a $1$--homological model for
a CDGA:
\begin{itemize}
\item The contraction defined in \cite{EM53,EM54} from $\bar{B}(A\ot A')$
to $\bar{B}(A)\ot \bar{B}(A')$, where $A$ and $A'$ are two
CDGAs.

$$C_{\scst \bar{B}\ot}:\{\bar{B}(A \ot A'),\bar{B}(A)
\ot \bar{B}(A'), f_{\scst \bar{B}\ot}, g_{\scst \bar{B}\ot},
\phi_{\scst \bar{B}\ot}\}\,;$$

\begin{itemize}
\item $f_{\scst \bar{B}\ot}[a_1 \otimes a'_1| \cdots |a_n \otimes
a'_n]$
\begin{eqnarray*}
\begin{array}{cl}
&=\displaystyle\sum _{i=0}^n
\xi_{\scst A} (a_{\scst i+1} *_{\scst A}\cdots a_{\scst n})
\xi_{\scst A'}(a'_{\scst 1}*_{\scst A'}
\cdots a'_{\scst i})[a_{\scst 1}| \cdots |a_{\scst i}]
\otimes [a'_{\scst i+1}| \cdots
|a'_{\scst n}]
\end{array}\end{eqnarray*}

\item $g_{\scst \bar{B}\ot}
([a_1| \cdots |a_n] \ot [a'_1| \cdots |a'_m] )$
\begin{eqnarray*}
=[a_1 \otimes \theta '| \cdots |a_n \otimes \theta
']\star [ \theta \otimes a'_1| \cdots | \theta \otimes
a'_n]\,,
 \end{eqnarray*}

where $\theta$ and $\theta '$ are the units in $A$ and $A'$ respectively.\\

\item up to sign,
$\phi_{\scst \bar{B}\ot}([a_1\ot a'_1|\cdots|a_n\ot a'_n])$
\begin{eqnarray*}
\begin{array}{cll}
\hspace{0.2cm}&=\displaystyle\sum&
\pm \xi_{\scst A} (a_{\scst n-q+1}*_{\scst A}\cdots  a_{\scst n})
[a_{\scst 1}\ot a'_{\scst 1}|\cdots |
a_{\scst \bar{n}-1}\ot a'_{\scst \bar{n}-1}\\&&
|a'_{\scst \bar{n}}*_{\scst A'}\cdots a'_{\scst n-q}|
c_{\scst \pi (0)}|\cdots |c_{\scst  \pi (p+q)}]\,,
\end{array} \end{eqnarray*}

\noindent where $\bar{n}=n-p-q$,
$(c_{\scst 0}, \ldots ,c_{\scst p+q})=(a_{\scst \bar{n}},
 \ldots ,a_{\scst n-q}, a'_{\scst n-q+1} ,\ldots a'_{\scst n})$
 and the sum is taken over all the
 $(p+1,q)$--shuffles $\pi$ and $0\leq p\leq n-q-1\leq n-1$.\\

Let us note that the complexity of $g_{\scst \bar{B}\ot}$
 and $\phi_{\scst \bar{B}\ot}$ is exponential since shuffles are involved
 in both formulas.
\end{itemize}

Given a tensor product $\ot _{i\in I}A_i$ of CDGAs, a
contraction from $\bar{B}(\ot_{i\in I}A_i)$ to $\ot_{i\in I}\bar{B}(A_i)$ is
easily determined by applying $C_{\scst \bar{B}\ot}$ several times in
a suitable way. This new contraction is also denoted by $C_{\scst
\bar{B}\ot}$.

\item The isomorphism of differential graded algebras (therefore, a contraction)
$$C_{\scst \bar{B}E}:
\{\bar{B}(E(u,2n+1)),\Gamma(\underline{u},2n+2),
f_{\scst \bar{B}E},g_{\scst \bar{B}E},0\} $$
described in \cite{EM54}, where
$$f_{\scst \bar{B} E}
([u|\stackrel{\mbox{\scriptsize $m$ times}}\cdots |u])
=\underline{u}^{(m)};\,\,\,
 g_{\scst \bar{B} E}(\underline{u}^{(m)})
 =[u|\stackrel{\mbox{\scriptsize $m$ times}}\cdots |u]\,.$$

\item The contraction
$$C_{\scst \bar{B}P}:\{\bar{B}(P(v,2n)),E(\underline{v},2n+1),
f_{\scst \bar{B}P},g_{\scst \bar{B}P},
\phi_{\scst \bar{B}P}
\}$$
stated in \cite{EM54}, where
$$f_{\scst \bar{B}P}([v^r])=\left\{\begin{array}{ll}
0&\mbox{if }r\neq 1\\
\underline{v}&\mbox{if } r=1
\end{array}\right.,\,\,\,
 f_{\scst \bar{B}P}([v^{\scst r_1}|\cdots |v^{\scst
r_m}])=0\,;$$
$$
g_{\scst \bar{B}P}(\underline{v})=[v]\,\,\, \mbox{ and } \,\,
 \phi_{\scst \bar{B}P}([v^{\scst
r_1}|\cdots |v^{\scst r_m}])=[v|v^{\mbox{$\scst r_1-1$}}|\cdots |v^{\scst
r_m}]\,.$$
\end{itemize}

Thanks to these three contractions, it is possible to establish, by
composition and tensor product of contractions, the
following semifull algebra contraction
$C=(f,
g, \phi)$:
$$\bar{B}(\ot_{i\in I} A_i) \; \Rightarrow \; \ot_{i\in I} \bar{B}(A_i) \;
\Rightarrow \; \ot_{i\in I} hBA_i\,,$$
where $hBA_i$ represents an exterior or a divided power algebra with
a generator $\underline{x}_i$,
depending on whether $A_i$ is a polynomial or an exterior algebra with
a generator $x_i$.

 In order to obtain the differential
structure of the $1$--homological model
for the twisted tensor product $\tilde{\ot}^{\rho}_{i\in I} A_i$,
the next step is to perturb $C$.
 The perturbation $\rho$
produces a perturbation--derivation $\delta$ on the tensor
differential of $\bar{B}(\ot_{i \in I} A_i)$:
$$\displaystyle\delta([a_1|\cdots|a_n])=\sum_{i=1}^n
(-1)^{|[a_1|\cdots|a_{i-1}]|}[a_1|\cdots|\rho(a_i)|\cdots|a_n]\,.$$

Now, by applying Theorem \ref{lph}, a new semifull algebra
contraction $(f_{\delta},g_{\delta},\phi_{\delta})$ is constructed:

$$\bar{B}(\tilde{\ot}^{\rho}_{i\in I} A_i)
\stackrel{(C)_{\delta}} {\Rightarrow} (\ot_{i \in I}
hBA_i,d_{\delta})\,,$$

\noindent where the differential $d_{\delta}$ is determined by the perturbation
procedure (Basic Perturbation Lemma). That means that
$hBA=\ot_{i \in I} (hBA_i, d_{\delta})$ is a
{\bf $1$--homological model} of $A=\tilde{\ot}_{i \in I} A_i$. Let
us emphasize that the Basic Perturbation Lemma provides finite formulas. Indeed, this is a
consequence of two facts: the perturbation $\delta$ does not
change the simplicial degree and $\phi$ increases this degree.

\begin{procedure}Algorithm {\bf Alg1}.

\vspace{0.5cm}

\hspace{-0.5cm}\begin{tabular}{ll}
\hline\\
{\bf Input}: & A finite CDGA $A$:
$((n_1,\dots,n_k),(v_1,\dots,v_k))$.\\\\
{\bf  Output:} & $((n_1+1,\dots,n_k+1),(w_1,\dots,w_k))$\\
&a $1$--homological model of
the CDGA $\ot_{i=1,\ldots, k}^{\rho} A_{i}$,
$A_i$ being the\\
& exterior algebra
$ E(x_i,n_i)$, if $n_i$ is odd and $P(x_i,n_i)$ if $n_i$ is
even.\\\\
\hline
\end{tabular}

\begin{tabbing}
$w_1=0$,\\
{\bf for} \= $i=2$ {\bf to} $k$\\
\>$w_i=d_{\delta}(\underline{x}_i)$, where $x_i$ is the algebra generator
of degree $n_i$\\
{\bf endfor}\\
\end{tabbing}
\end{procedure}

  Naturally, the first components of the vector $\bar{v}$ must be
  zero, because they correspond to the image of the
  algebra generators with the lowest degree under $\rho$.

Moreover, a general algorithm for computing the
$1$--homology of
CDGAs can be described . Clearly, the homology of the $1$--homological model obtained can be
computed using an algorithm based on the establishment of  Smith's
normal form of the matrices representing the differentials at each degree
\cite{Veb31,Mun84}.

The computational cost of constructing the contraction $(C)_{\delta}$ is
high. Let us note that both the inclusion and homotopy operators of the
contraction $C_{\scst \bar{B} \ot}$ give an answer in exponential time.
In fact, the formula of the differential operator
$d_{\delta}$ produced by the homological perturbation machinery is given
by:
\begin{eqnarray*}d_{\delta} \;\; = \;\;
f \, \delta  ( 1 - \phi\, \delta +
  \phi\,  \delta\, \phi\, \delta -
\cdots ) \, g\,. \end{eqnarray*}

 With regard to the previous remarks,
a first impression is that obtaining
$d_{\delta}$ generally becomes  a procedure
of exponential nature.

It is possible to take advantage of $d_{\delta}$
being a derivation. Indeed, the fact that $d_{\delta}$ is a derivation implies
that it is only necessary to know this  morphism
applied to the generators of the model (let us observe that there are
as
many generators as the cardinal of the  set of indices $I$ indicated).
This is an enormous improvement in the computation of the
differential on the small model. In spite of this,  computing  $d_{\delta}$ on an
algebra generator is  extremely time--consuming.

\section{Inversion Theory}\label{inversiones}

In this section, we  go
further in the simplification of the
computation of the differential $d_{\delta}$.
For clarity, we begin this work considering only two algebras.

As we have seen before, obtaining $d_{\delta}$ is an extremely expensive procedure.
  The morphism responsible for this
is the homotopy
operator, $\phi$, due, essentially, to the shuffles that are
involved in the formulas of $\phi_{\scst B\ot}$ and $g_{\scst B\ot}$. We
intend to eliminate these shuffles, and, with this aim in mind, we define
the concept of inversion.

\begin{definition}{\em
Let $A$ and $A'$ be CDGAs and let us consider a  homogeneous
element
$[a_1\ot a'_1|a_2\ot a'_2|\cdots |a_n\ot a'_n]$ from $\bar{B}(A\ot
A')$. We say that a component $\theta\ot a'_i$ from that element,
is responsible for an {\em inversion}, if
 there exists an index  $j> i$ with $a_j\neq
\theta  $ (where $\theta$ is the unit of $A$). In this sense,
such an element presents $k$ inversions if there exist $k$
components responsible for an inversion.
 }\end{definition}

We will say that an element from $\bar{B}(A\ot
A')$ has $k $ inversions, if it is a sum of elements which each
have, as a minimum,
$k$ inversions.

\vspace{0.3cm}

Let us consider the contraction
$$(f_{\scst \bar{B}\ot},g_{\scst \bar{B}\ot},\phi_{\scst \bar{B}\ot} ):
\bar{B}(A\ot A')\Rightarrow \bar{B}(A)\ot \bar{B}(A')$$
described in the previous section.
We analyze
the behaviour of the component morphisms of this contraction with respect to
 inversions. For this purpose, we do not take into account the signs
in the formulas referred to.

\begin{itemize}
\item The image of an element
with at least one inversion under $f_{\scst \bar{B}\ot}$,  is null.

\item The injection $g_{\scst \bar{B}\ot }$,  applied to
$[a_1|\cdots |a_n]\ot [a'_1|\cdots|a'_m]$,
produces:
\begin{itemize}
\item a unique term with no inversions (that one which comes from
juxtaposition),
\item $n$ terms with one inversion,
\item $\left(\begin{array}{c} n+m \\
n\end{array} \right)-n-1$ terms with more than one inversion.
\end{itemize}

\item As for the homotopy operator $\phi_{\scst \bar{B}\ot}$,
we can state that
the image  of  a homogenous
element  under $\phi_{\scst \bar{B}\ot}$ gives rise to a sum of elements
which, if non null, have at least one more inversion than the
original one.
Let us note that  an inversion is
produced by the component
$\;a'_{\scst \bar {n}}*_{\scst A'}\cdots a'_{\scst n-q}$,
which is always on the left side of those components $a_{\scst \bar{n}}, \ldots ,a_{\scst
n-q}$ of each summand in the
formula of $\phi_{\scst \bar{B}\ot}$.
\end{itemize}
\vspace{0.5cm}

Let us consider the contraction which provides us with a $1$--homological
model for the tensor product of two CDGAs, $A$ and $A'$:

\begin{equation}\label{ctotal}
(f,g,\phi): \bar{B}(A\ot A')\Rightarrow \bar{B}(A)\ot\bar{B}(A')
\Rightarrow hBA\ot hBA'
\end{equation}
where $$\begin{array}{l}
f=(f_{\scst \bar{B}A}\ot f_{\scst \bar{B}A'})f_{\scst
\bar{B}\ot},\\
g=g_{\scst \bar{B}\ot}(g_{\scst \bar{B}A}\ot g_{\scst
\bar{B}A'})\\
\phi=\phi_{\scst \bar{B}\ot} +
g_{\scst \bar{B}\ot}(\phi_{\scst \bar{B}A}\ot g_{\scst \bar{B}A'}f_{\scst \bar{B}A'}+
1_{\scst \bar{B}A}\ot\phi_{\scst \bar{B}A'})f_{\scst \bar{B}\ot}
\end{array}$$

Let us note that the image of  an element with an
inversion under $f$ is also null, since the first morphism applied is $f_{\scst
\bar{B}\ot}$.

\vspace{0.3cm}

 Now we assume that there is a perturbation $\rho$ of the tensor product
of the algebras $A$ and $A'$. This perturbation induces, in a natural way,
a perturbation $\delta$ on
$\bar{B}(A\ot A')$. Let us  analyze the behaviour of such a
morphism with respect to inversions.

\begin{lemma}\label{deltadeinv}
Let us consider a perturbation $\delta$ for $\bar{B}(A\ot A')$
induced by a perturbation--derivation $\rho$ for $A\ot A'$ such that
$\rho (A)\subset A$. The
 image of a homogeneous element
 with $k$ inversions under $\delta$, is a sum of elements with at least
$k-1$ inversions.
\end{lemma}

\noindent{\bf Proof.}

Let us point out that a component of a
homogeneous element from $\bar{B}(A\ot A')$ is responsible for, at
most, one inversion and that $\delta$ acts only on a component of the
element at each term of the resultant sum.

\hfill{$\Box$}

Attending to the Basic Perturbation Lemma, one can obtain from the
contraction (\ref{ctotal}), a new contraction:

\begin{eqnarray*}
(f_{\delta},g_{\delta},\phi_{\delta}):
 \bar{B}(A\tilde{\ot}^{\scst \rho} A')
\Rightarrow  (hBA\ot hBA',d_{\delta})\,.
\end{eqnarray*}

We recall the formula for $d_{\delta}$:

$$d_{\delta} \;\; = \;\;
f \, \delta \, ( 1 - \phi\, \delta +
  \phi \,\delta\, \phi\, \delta -
\cdots ) \, g.  $$

We can observe that $f$ is the last morphism applied in the formula.
If at any stage,
an element $y$ obtained by applying $\phi$,
has more than one inversion, then $\delta (y)$
will have at least one inversion.
In this way, each time we apply $\delta\,\phi$, we
obtain a sum of homogeneous elements with at least one inversion, and, therefore, the
image of these elements under  $f$ is null. This means that we only have to
consider the summands of $\phi$ having, at most, one
inversion.

In consequence, we can establish the following theorem
where we considerably reduce
the complexity of the computation of $d_{\delta}$.

\begin{theorem}\label{reduccion}
The formula for $\phi$,
that is involved in the definition of $d_{\delta}$,
is the following:

$$\phi=\bar{\phi}_{\scst \bar{B}\ot}+\bar{g}_{\scst \bar{B}\ot}
(\phi_{\scst \bar{B}A }\ot
g_{\scst \bar{B}A'}f_{\scst \bar{B}A'}+1\ot \phi_{\scst \bar{B}A'})
f_{\scst \bar{B}\ot},
$$
where
\begin{itemize}
\item
$\bar{\phi}_{\scst \bar{B}\ot}([a_1\ot a'_1|\cdots|a_n\ot a'_n])$\\\\
$\begin{array}{cll}
=&\displaystyle\sum_{\scst 0\leq p\leq n-q-1\leq n-1}&
(-1)^{\varphi(n,p,q)}\xi_{\scst A} (a_{\scst n-q+1}*_{\scst A}\cdots a_{\scst n})
[a_{\scst 1}\ot a'_{\scst 1}|\cdots |
a_{\scst \bar{n}-1}\ot a'_{\scst \bar{n}-1}\\&&
|a'_{\scst \bar{n}}*_{\scst A'}\cdots a'_{\scst n-q}|
a_{\scst\bar{n}}|\cdots |a_{\scst n-q}|a'_{\scst n-q+1}|\cdots |a'_{\scst n}]
\end{array}$\\\\

\noindent being $\bar{n}=n-p-q$ and
\begin{eqnarray*}
\varphi (n,p,q)&=&\bar{n}-1+
|[a_1|\cdots |a_{\bar{n}-1}]|_t+|[a'_1 |\cdots |a'_{n-q}]|_t\\
&&+\sum_{k=0}^p\sum_{\ell =0}^k|a_{n-q-k}|\;|a'_{n-q-\ell }|\,.
\end{eqnarray*}

\item  $\bar{g}_{\scst
\bar{B}\ot}([a_1|\cdots|a_n]\ot[a'_1|\cdots|a'_m])$\\\\
$\begin{array}{cl}
\hspace{1cm}=&[a_1|\cdots|a_n|a'_1|\cdots|a'_m]\\\\
&+\displaystyle\sum_{\scst i=1}^{\scst n-1}(-1)^{|[a_{i+1}|\cdots
|a_n]||a'_1|}\;
[a_1|\cdots|a_i|a'_1|a_{i+1}|\cdots |a_n|a'_2|\cdots |a'_m]\\\\
&+(-1)^{|a'_1||[a_1|\cdots|a_n]|}\;
[a'_1|a_1|\cdots|a_n|a'_2|\cdots|a'_m]\,.
\end{array}$
\end{itemize}
\end{theorem}

Let us note that now the number of summands in
the formula above for $\phi_{\scst \bar{B}\ot}$ is
$$\sum_{q=0}^{n-1}\sum_{p=0}^{n-q-1}1=\frac{n^2+n}{2}\,,$$ in contrast to
the original
number of summands: $$\sum_{q=0}^{n-1}\sum_{p=0}^{n-q-1}\left(\begin{array}{c}
p+q+1\\q
\end{array}\right)\;\; = \;\; 2^{n+1}-n-2\,.$$
On the other hand, the formula for $g_{\scst \bar{B}\ot}$ is
reduced to $n$ summands, instead of $\left(\begin{array}{c}
m+n\\n
\end{array}\right)$.

This theorem is easy to generalize, by induction, to the general case
of a twisted tensor product of CDGAs  $\tilde{\ot}_{i\in I}^{\rho}A_i$ with
$I=\{1,\ldots ,n\}$,
where $A_i$ is an exterior or a polynomial algebra with generator $x_i$
and $|x_i|\leq |x_{i+1}|$. Therefore, as we
saw in Section \ref{algoritmo},
the following semifull
algebra contraction can be established:
\begin{eqnarray*}
 \bar{B}(\ot_{i=1}^n A_i )\Rightarrow
\ot_{i=1}^n\bar{B}(A_i)\Rightarrow \ot_{i=1}^n hBA_i
\end{eqnarray*}
where $hBA_i$ is a polynomial or a divided power algebra.

The key to understanding the generalization is the fact that
the inversions
in $\bar{B}(\ot_{i=1}^{n}A_i)$ are those of the last tensor
product (as they were defined at the beginning of the Section,
with $A=\ot_{i=1}^{n-1}A_i$ and $A'=A_n$) along with those of
$\ot_{i=1}^{n-1}A_i=(\ot_{i=1}^{n-2}A_i)\ot A_{n-1}$ with respect to
the last tensor product, and so on.

  Summing up, we obtain an algorithm {\bf Alg2} having the
  same input and output as the Algorithm {\bf Alg1} of the last
  section, but speeding up the steps concerning the image of the algebra
  generators  under $d_{\delta}$.

\begin{procedure}Algorithm {\bf Alg2}.

\vspace{0.5cm}

\hspace{-0.5cm}\begin{tabular}{c}
\hline\\
{\bf Input} and {\bf Output}: the same as in {\bf Alg1}.\\\\
\hline
\end{tabular}

\begin{tabbing}
$w_1=0$,\\
{\bf for} \= $i=2$ {\bf to} $k$\\
\>$w_i=d_{\delta}(\underline{x}_i)$, where $x_i$ is the algebra generator
of degree $n_i$ \\
\>(using Theorem 4)\\
{\bf endfor}\\
\end{tabbing}
\end{procedure}

\section{Complexity}\label{complejidad}

  In this section we give a comparison of the algorithms {\bf
  Alg1}
and {\bf Alg2} from the point of view of their complexity. We are mainly
interested in measuring the efficiency of the corresponding steps
concerning the obtention of the differential $d_{\delta}$.
We consider the degree of the algebra generator as
 the size of an instance. We take as elementary
operations those ones generating each homogeneous term produced by the
different morphisms and a worst--case analysis of the algorithms is carried out.

  We calculate the total number
  of elementary operations needed for
  computing $d_{\delta}$ on a generator $\underline{x}_{k}$
   of degree
  $k$, for both {\bf Alg1} and {\bf Alg2}. We hold that this number
  for {\bf Alg1} is
  $$\displaystyle\sum_{i=0}^{\left\lfloor k/k_0+1\right\rfloor}
  (i!\,r^i+((i+1)!\,+(i+2)!\,)r^{i+1})
  \left(\prod_{j=1}^i 2^{j+1} -j-2\right)\, ,$$
and for {\bf Alg2} is
$$\displaystyle\sum_{i=0}^{\left\lfloor k/k_0+1\right\rfloor}
(r^i + (i+3)r^{i+1})
  \left(\prod_{j=1}^i \frac{j^2+ j}{2}\right)\, ,$$
where $k_0=\displaystyle \mbox{min}_{\scst 1\leq i\leq n}|x_i|$
 and $r$ is the maximum number of summands given by
$\rho(x_i)$, where $x_i$ ranges over all the algebra generators
of $A$.

    In the following table, the required time for
    computing $d_{\delta}(\underline{x}_k)$ is showed, supposing that our
    computer carries out $10^6$ elementary operations per second. Let us denote
     $s = \left\lfloor k/k_0+1\right\rfloor$. Note that,
     for example, that $s=5$ means that if $k_0=10$, the
degrees of the algebra generators range over the set
$\{10,11,\ldots ,54,55\}$.

  \begin{table}\caption{Required time}
\begin{center}
\begin{tabular}{lclclcl}
\hline\noalign{\smallskip}
&\hspace{0.5cm}& &\hspace{1cm}&in {\bf Alg1}& \hspace{1cm}& in {\bf Alg2} \\
\noalign{\smallskip}\hline\noalign{\smallskip}
$s=3;$&&$ r=2$ & &0.1 sec. && 0.002 sec. \\
$s=3;$&&$ r=3$&& 0.5 sec. && 0.009 sec.\\
$s=4;$&&$r=2$&& 31.29 sec. & &0.04 sec. \\
$s=5;$&&$ r=3$&& 3.19 days && 17.6 sec.\\
$s=6;$&&$ r=2$&& 1.45 years&& 1.17 min \\
$s=6;$&&$ r=3$& &24.75 years & &19.56 min\\
\hline\end{tabular}\end{center}\end{table}

\section{Implementation Performance}\label{ejemplos}

    The algorithm {\bf Alg2} has been  implemented. The user supplies an
    encoding of a finite CDGA in the form of
    a twisted tensor product of exterior and polynomial algebras to the
    program, which computes an encoding of a small
    $1$--homological model of this algebra.

  This program is written in Mathematica 3.0, consisting in $300$
     lines of code and  $10$ basic functions.

 In order to give some indication of the
 implementation, we report on the time taken to compute $1$--homological
 models for certain CDGAs:

\begin{enumerate}
\item $E(x_1,x_2,x_3;1)\ot P(x_4,x_5;2)\ot
E(x_6;3)\ot P(x_7; 4)\ot P(x_8;6)$.
\newline\begin{tabbing}
{\bf Input}: \=
$((1,1,1,2,2,3,4,6),$\\
\>$(0,0,0,x_{\scst 1}-x_{\scst 2},x_{\scst 2},
x_{\scst 1}x_{\scst 2},x_{\scst 2}x_{\scst 4}+x_{\scst 1}x_{\scst 5}+
x_{\scst 1}x_{\scst 2}x_{\scst 3},
x_{\scst 1}x_{\scst 2}x_{\scst 6}))$.\\\\
{\bf  Output:} \= $((2,2,2,3,3,4,5,7),(0,0,0,\underline{x}_{\scst 1}-
\underline{x}_{\scst 2}
,\underline{x}_{\scst 2},0,\underline{x}_{\scst 1}\underline{x}_{\scst 2}
-2\underline{x}_{\scst 2}^2,0))$\\\\
{\bf Time: }\= $d_{\delta}(\underline{x}_4)$ in 0.55 sec.\\
\> $d_{\delta}(\underline{x}_5)$ in 0.27 sec.\\
\> $d_{\delta}(\underline{x}_6)$ in 0.33 sec.\\
\> $d_{\delta}(\underline{x}_7)$ in 2.14 sec.\\
\> $d_{\delta}(\underline{x}_8)$ in 0.28 sec.\\
\end{tabbing}
\item $E(x_1;1)\ot P(x_2;2)\ot P(x_3;6)\ot  P(x_4; 10)
\ot P(x_5;26)$.
\newline\begin{tabbing}
{\bf Input}: \=
$((1,2,6,10,26),(0,-2x_{\scst 1},x_{\scst 1}x_{\scst 2}^2,
3x_{\scst 1}x_{\scst 2}x_{\scst 3},8x_{\scst 1}x_{\scst 2}x_{\scst 3}^2
x_{\scst 4}))$.\\\\
{\bf  Output:} \= $((2,3,7,11,27),(0,-2\underline{x}_{\scst 1},-8
\underline{x}_{\scst 1}^3
,-192\underline{x}_{\scst 1}^5,21799895040\underline{x}_{\scst 1}^{13}))$\\\\
{\bf Time: }\=$d_{\delta}(\underline{x}_2)$ in 0.16 sec.\\
\> $d_{\delta}(\underline{x}_3)$ in 0.60 sec.\\
\> $d_{\delta}(\underline{x}_4)$ in 2.36 sec.\\
\> $d_{\delta}(\underline{x}_5)$ in 1571.47 sec.\\
\end{tabbing}
\item $E(x_1;1)\ot P(x_2;2)\ot P(x_3;4)\ot
E(x_4;5)\ot E(x_5; 7)\ot P(x_6;14)$.
\newline\begin{tabbing}
{\bf Input}: \=
$((1,2,4,5,7,14),(0,-x_{\scst 1},x_{\scst 1}x_{\scst 2},0,
2x_{\scst 1}x_{\scst 4},-x_{\scst 1}x_{\scst 4}x_{\scst 5}))$.\\\\
{\bf  Output:} \= $((2,3,5,6,8,15),(0,-\underline{x}_{\scst 1},
-\underline{x}_{\scst 1}^2,0,0))$\\\\
{\bf Time: }\= $d_{\delta}(\underline{x}_2)$ in 0.17 sec.\\
\> $d_{\delta}(\underline{x}_3)$ in 0.33 sec.\\
\> $d_{\delta}(\underline{x}_4)$ in 0.26 sec.\\
\> $d_{\delta}(\underline{x}_5)$ in 0.17 sec.\\
\end{tabbing}
\end{enumerate}

For each algebra, we summarize the results and the
 time taken to compute a description of the model. All CPU times
 are in seconds and calculations were carried out on a Pentium III,
 128Mb RAM, 7.2Gb Hard disk space.

  This program produces as output an encoding of a
   certain differential graded algebra which could be introduced
   into another program in order to calculate the homology of such
   objects. We intend to tackle this task  in the near future.

\end{document}